\documentstyle[11pt]{article}

\title{NONSTANDARD TRANSFINITE GRAPHS AND NETWORKS OF HIGHER RANKS}
\author{A. H. Zemanian}
\date{}

\begin{document}
\newcommand{\N} {I \kern -4.5pt N}
\newcommand{\R} {I \kern -4.5pt R}
\maketitle
\baselineskip21pt

{\ Abstract --- In Chapter 8 of the book, ``Graphs and Networks: Transfinite 
and Nonstandard,'' (published by Birkhauser-Boston in 2004), 
nonstandard versions of transfinite graphs 
and of electrical networks having 
such graphs were defined and examined but only for the first two ranks,
0 and 1, of transfiniteness. 
In the present work, these results are extended to 
higher ranks of transfiniteness.  Such is done in detail for the 
natural-number ranks and also for the first transfinite ordinal rank $\omega$. 
Results for still higher ranks of transfiniteness can be 
established in much the same way.  Once the transfinite 
graphs of higher ranks are established, theorems concerning 
the existence of hyperreal operating points and the satisfaction 
of Kirchhoff's laws in nonstandard networks of higher ranks
can be proven just as they are 
for nonstandard networks of the first rank. \\

} 

\section{Introduction}

In a prior publication \cite[Chapter 8]{gn} we defined and 
examined nonstandard versions of graphs that are conventionally infinite 
as well as those that are transfinite but only of the first rank of 
transfiniteness. We also examined nonstandard, resistive, electrical 
networks having such graphs and established an existence theorem for 
their operating points (i.e., their hyperreal current-voltage regimes)
as well as Kirchhoff's laws for those nonstandard networks.
In this work, we shall extend these results to graphs and networks 
having higher ranks of transfiniteness.  We do so in detail 
for the natural-number ranks and also 
for the first transfinite-ordinal rank $\omega$.
Results for still higher ranks of transfiniteness can be established in 
virtually the same way; the development for the 
successor-ordinal ranks (resp. limit-ordinal ranks) are 
virtually the same as that for the natural-number
ranks (resp. the rank $\omega$).

All this is accomplished through a recursive analysis 
proceeding along increasing ordinal ranks.  The first two steps of that
recursion concern the ranks 0 and 1.  These have been
explicated in \cite[Chapter 8]{gn} and will not be repeated here.  
Our notation and terminology is the same as that used
in \cite{gn}.

\section{Nonstandard $\mu$-Graphs}

Let $\mu$ be a natural number no less than 2.  Our development of a 
nonstandard $\mu$-graph starts with a given sequence
$\langle G_{n}^{\mu}\!: n\in \N\rangle$, where 
\[ G_{n}^{\mu}\;=\;\{X_{n}^{0},B_{n},X_{n}^{1}.\ldots,X_{n}^{\mu}\} \]
is a standard transfinite graph of rank $\mu$.
Here, we are defining the branches (i.e., the members of $B_{n}$)
as pairs of 0-nodes (i.e., members of $X_{n}^{0}$).
This differs from the definition of branches given in \cite[page 6]{gn}
based upon elementary tips but only in a nonessential way.
We can indeed use all the ideas and results given in 
\cite[Chapter 2]{gn}.\footnote{In this regard, see how the definition 
$G^{1}=\{X^{0},B,X^{1}\}$ of a 1-graph given in \cite[page 163]{gn}
differs from the definition 
${\cal G}^{1}=\{{\cal B},{\cal X}^{0},{\cal X}^{1}\}$ of a 1-graph given in 
\cite[page 8]{gn} in that nonessential way.}  Thus,
\[ G_{n}^{\mu-1}\;=\;\{X_{n}^{0},B_{n},X_{n}^{1},\ldots,X_{n}^{\mu-1}\} \]
is the $(\mu-1)$-graph of $G_{n}^{\mu}$.

The {\em extremities} of $G_{n}^{\mu-1}$ are taken to be its 
$(\mu-1)$-tips and also the exceptional elements of the $\mu$-nodes
of $G_{n}^{\mu}$.  (The exceptional element, if it exists, 
of a $\mu$-node $x^{\mu}$ is the unique node of rank less than 
$\mu$ contained in $x^{\mu}$;  see 
\cite[page 11]{gn}.)  Let ${\cal T}_{n}^{\mu-1}$ be the set of  
$(\mu-1)$-tips of $G_{n}^{\mu-1}$.  Let a typical $\mu$-node of 
$G_{n}^{\mu}$ be denoted by $x_{n,k}^{\mu}$.  We are 
indexing those $\mu$-nodes 
by $k$, and we let $K$ be the index set for those $\mu$-nodes.
In accordance with the 
partitioning defined by the $x_{n,k}^{\mu}$ , ${\cal T}_{n}^{\mu-1}$
is partitioned into subsets ${\cal T}_{n,k}^{\mu-1}$, where 
${\cal T}_{n,k}^{\mu-1}$ is the set of all the $(\mu-1)$-tips in 
$x_{n,k}^{\mu}$.  Thus, 
\[ {\cal T}_{n}^{\mu-1}\;=\;\cup_{k\in K}{\cal T}_{n,k}^{\mu-1}, \] 
where $K$ serves also as the index set for that partitioning.

If a $\mu$-node $x_{n,k}^{\mu}$ of $G_{n}^{\mu}$ has an 
exceptional element $x_{n,k}^{{\alpha}}$ $({\alpha} < {\mu})$,
let ${\cal Z}_{n,k}$ denote the singleton set
${\cal Z}_{n,k}=\{x_{n,k}^{{\alpha}}\}$.  Otherwise, let
${\cal Z}_{n,k}\;=\;\emptyset$.  In either case, by the definition of any 
standard $\mu$-node $x_{n,k}^{\mu}$, 
${\cal Z}_{n,k} \cap{\cal Z}_{n,l} =\emptyset$ whenever $k\neq l$,
and we have 
\[ x_{n,k}^{\mu}\;=\;{\cal T}_{n,k}^{\mu-1} \cup {\cal Z}_{n,k}. \]

If $e_{n}$ and $f_{n}$ are two extremities in the same $\mu$-node 
$x_{n,k}^{{\mu}}$ of $G_{n}^{{\mu}}$, we say that $e_{n}$ and $f_{n}$ are 
{\em shorted together}, and we write $e_{n} \asymp f_{n}$ to 
denote this fact.  

Out next objective is to make an ultrapower construction of 
the nonstandard $\mu$-nodes and thereby obtain the nonstandard
${\mu}$-graph $^{*}\!G^{{\mu}}$.  We already have at hand the nonstandard 
0-graph $^{*}\!G^{0}=\{\,^{*}\!X^{0},\,^{*}\!B\}$ 
\cite[page 155]{gn} and the nonstandard 1-graph $^{*}\!G^{1}=
\{\,^{*}\!X^{0},\,^{*}\!B,\,^{*}\!X^{1}\}$ \cite[page 164]{gn}.
So recursion will yield the nonstandard $\mu$-graphs
\begin{equation}
^{*}\!G^{{\mu}}\;=\;\{\,^{*}\!X^{0},\,^{*}\!B,\,^{*}\!X^{1},\ldots,\,^{*}\!X^{\mu}\},  \label{2.1}
\end{equation}
where $^{*}\!X^{{\mu}}$ is the set of nonstandard $\mu$-nodes.

To this end, let $\cal F$ be any chosen and fixed nonprincipal ultrafilter.
Let $\langle e_{n}\rangle$ be a sequence where each $e_{n}$ is 
an extremity of $G_{n}^{{\mu}-1}$.  Two such sequences 
$\langle e_{n}\rangle$ and $\langle f_{n}\rangle$ are said 
to be {\em equivalent} if $e_{n}=f_{n}$
for almost all $n$ (modulo $\cal F$); i.e., $\{n\!: e_{n}=f_{n}\}\in{\cal F}$.
This partitions the set of all extremities into equivalence classes;
indeed, reflexivity and symmetry are obvious and for transitivity
we have, with $\langle g_{n}\rangle$ being another sequence of 
extremities, 
\[ \{n\!: e_{n}=f_{n}\}\,\cap\,\{n\!: f_{n}=g_{n}\}\,\subseteq\{n\!: f_{n}=g_{n}\} \]
so that $\langle e_{n}\rangle$ is also equivalent to $\langle g_{n}\rangle$.
Each such equivalence class will be called a {\em nonstandard extremity} and 
denoted by ${\bf e}=[e_{n}]$, where $\langle e_{n}\rangle$ is any 
representative of that equivalence class.

Given any sequence $\langle e_{n}\rangle$ of extremities,
let $N_{t^{{\mu}-1}}$ be the set of all $n$ for which $e_{n}$ is a 
$({\mu}-1)$-tip of $G_{n}^{{\mu}-1}$, and let $N_{x}$ be the set of all 
$n$ for which $e_{n}$ is an exceptional element of a ${\mu}$-node of
$G_{n}^{\mu}$. Consequently, $N_{t^{\mu-1}}\cup N_{x}=\N$ and
$N_{t^{\mu-1}}\cap N_{x}=\emptyset$.  So, exactly one of 
$N_{t^{\mu-1}}$ and $N_{x}$ is a member of $\cal F$.
If it is $N_{t^{\mu-1}}$ (resp. $N_{x}$), we define $\langle e_{n}\rangle$
as being a representative of a {\em nonstandard $({\mu}-1)$-tip}
${\bf t}^{{\mu}-1}=[e_{n}]$
(resp. a representative of a {\em nonstandard exceptional element}) ${\bf x}
=[e_{n}]$).  

In the latter case of an exceptional element, the $e_{n}$ are nodes of 
$G_{n}^{{\mu}-1}$ for almost all $n$, but they need not be of the same rank;  
their ranks can vary through values no larger than ${\mu}-1$. 
There are no more than finitely many such ranks.  
Let $K$ be the finite set of such ranks, and let $F_{k}$ 
denote the set of all $n$ for which the rank has the value $k$.
The sets $F_{k}$ are finitely many, pairwise disjoint, and their union 
is a member of $\cal F$.  Therefore, exactly one of 
those sets $F(\rho)$ is a member of $\cal F$
\cite[page 19, fact (4)]{gn}.  Consequently, we can identify the rank of 
$\bf x$ as the rank $\rho$ of that unique set $F({\rho})$,
and so we may denote $\bf x$ as ${\bf x}^{\rho}$.

Next step:  Let ${\bf e}=[e_{n}]$ and ${\bf f}=[f_{n}]$ be two
nonstandard extremities.  Let $N_{ef}=\{n\!: e_{n}\asymp f_{n}\}$
and $N_{ef}^{c}=\{n\!: e_{n}\not\asymp f_{n}\}$.
So, exactly one of $N_{ef}$ and $N_{ef}^{c}$ is a member of ${\cal F}$.
If it is $N_{ef}$ (resp. $N_{ef}^{c}$), we say that $\bf e$ is 
{\em shorted to} $\bf f$, and we write ${\bf e}\asymp{\bf f}$
(resp. we say that $\bf e$ is {\em not shorted to} $\bf f$, and we write
${\bf e}\not\asymp {\bf f}$). Also, we take it that $\bf e$ is 
{\em shorted to itself}.  This shorting is an equivalence relation
for the set of all nonstandard extremities.  Indeed, reflexivity and
symmetry are again obvious, and transitivity follows from 
\[ \{ n\!: e_{n}\asymp f_{n}\}\,\cap\,\{n\!: f_{n}\asymp g_{n}\}
\,\subseteq \,\{n\!: e_{n}\asymp g_{n}\}. \]
The resulting equivalence classes are defined to be the 
{\em nonstandard ${\mu}$-nodes}, and we use the boldface notation 
${\bf x}^{\mu}$ to denote a typical one.

Various properties of standard nodes transfer directly to nonstandard nodes.
For instance, if the nonstandard node ${\mu}$-node ${\bf x}^{{\mu}}$
has a nonstandard exceptional element 
${\bf x}^{{\rho}}=[e_{n}]$ (${\rho}<{\mu}$),
we have that for almost all $n$, $e_{n}\asymp f_{n}$, where
${\bf f}=[f_{n}]$ is a nonstandard $({\mu}-1)$-tip in ${\bf x}^{\mu}$; 
that is, every nonstandard exceptional element is shorted to at
least one nonstandard $({\mu}-1)$-tip.
This also implies that every nonstandard ${\mu}$-node has at least one 
nonstandard $({\mu}-1)$-tip.

For similar reasons, the exceptional element of a nonstandard
${\mu}$-node cannot be the exceptional element of any other nonstandard 
${\mu}$-node, and no nonstandard ${\mu}$-node can have two or more nonstandard
exceptional elements.

Let $^{*}\!X^{{\mu}}$ denote the set of all nonstandard ${\mu}$-nodes as 
determined by the given sequence $\langle G_{n}^{{\mu}}\rangle$
of standard ${\mu}$-graphs.  By recursion we can now define the nonstandard 
${\mu}$-graph $^{*}\!G^{{\mu}}$ as the $({\mu}+2)$-tuple given by 
(\ref{2.1}) above.

\section{Nonstandard Graphs of Rank $\vec{\omega}$}

We now take it that our recursive construction of the 
$\mu$-graphs can be continued indefinitely through all the
natural-numbers ranks.  That is, given any sequence 
$\langle G_{n}^{\vec{\omega}}\rangle$ of standard 
$\vec{\omega}$-graphs,\footnote{See \cite[Section 2.5]{gn} for the definition
of such a $\vec{\omega}$-graph.} we can construct as in the preceding section 
each set $^{*}\!X^{\mu}$
of nonstandard $\mu$-nodes from the sequence 
$\langle G_{n}^{\mu}\rangle$, where $G_{n}^{\mu}$ is the 
$\mu$-graph of the $\vec{\omega}$-graph 
$G_{n}^{\vec{\omega}}$,
this being so for every $\mu\in\N$.
Furthermore, for the sake of simplicity, we shall assume that none of the 
$G_{n}^{\vec{\omega}}$ contains $\vec{\omega}$-nodes.  As a result, 
the nonstandard graph $^{*}\! G^{\vec{\omega}}$ 
we shall specify in a moment will not
possess any nonstandard $\vec{\omega}$-node.\footnote{It may be possible 
to construct an equivalence class $[x_{n}^{\vec{\omega}}]$ (modulo $\cal F$) 
of standard $\vec{\omega}$-nodes to obtain nonstandard $\vec{\omega}$-nodes 
${\bf x}^{\vec{\omega}}$, but there is a problem concerning the ranks 
of the embraced nodes of the $x_{n}^{\vec{\omega}}$.
We have not resolved this matter.}

Thus, we can now define the {\em nonstandard $\vec{\omega}$-graph}
$^{*}\! G^{\vec{\omega}}$ as the sequence
\begin{equation}
^{*}\!G^{\vec{\omega}}\;=\;\{^{*}\!X^{0},^{*}\!B,\ldots,^{*}\!X^{\mu},\ldots\},  \label{3.1}
\end{equation}
where the entries $^{*}\!X^{\mu}$ extend throughout all the 
natural-numbers $\mu\in\N$.

\section{Nonstandard Graphs of Rank $\omega$}

We now start with a given sequence $\langle G_{n}^{\omega}\rangle$ of
standard $\omega$-graphs:
\[ G_{n}^{\omega}\;=\;\{X_{n}^{0},B_{n},X_{n}^{1}\ldots,X_{n}^{\mu},\ldots,X_{n}^{\omega}\}, \]
where none of the $G_{n}^{\omega}$ has any $\vec{\omega}$-node---in accordance 
with our assumption in Section 3.
We have
\[ G_{n}^{\vec{\omega}}\;=\;\{X_{n}^{0},B_{n},X_{n}^{1}\ldots,X_{n}^{\mu},\ldots\} \]
as the $\vec{\omega}$-graph of $G_{n}^{\omega}$.  The {\em extremities} of
$G_{n}^{\vec{\omega}}$ are its $\vec{\omega}$-tips and the exceptional elements of the 
$\omega$-nodes of $G_{n}^{\omega}$.  Those exceptional elements, if they exist,
are nodes of $G_{n}^{\vec{\omega}}$ with natural-number ranks.  

Given any sequence $\langle e_{n}\rangle$ of extremities, one from 
each $G_{n}^{\vec{\omega}}$, let $N_{t^{\vec{\omega}}}$ be the set of all $n$ for which 
$e_{n}$ is an $\vec{\omega}$-tip, and 
let $N_{x}$ be the set of all $n$ for which 
$e_{n}$ is an exceptional element.  Thus, 
$N_{t^{\vec{\omega}}}\cap N_{x}=\emptyset$
and $N_{t^{\vec{\omega}}}\cup N_{x}=\N$.  Consequently, exactly one of 
$N_{t^{\vec{\omega}}}$ and $N_{x}$ is a member of $\cal F$.  If it is 
$N_{t^{\vec{\omega}}}$ (resp. $N_{x}$), 
$\langle e_{n}\rangle$ is a representative of a 
{\em nonstandard $\vec{\omega}$-tip} ${\bf t}^{\vec{\omega}}=[e_{n}]$ 
(resp. a representative
of a {\em nonstandard exceptional element} ${\bf x}=[e_{n}]$).
In either case, we also refer to ${\bf e}=[e_{n}]$ as a 
{\em nonstandard extremity}.

Now, let ${\bf e}=[e_{n}]$ and ${\bf f}=[f_{n}]$ be two nonstandard 
extremities.  Let $N_{ef}=\{n\!: e_{n}\asymp f_{n}\}$ and
$N_{ef}^{c}=\{n\!: e_{n}\not\asymp f_{n}\}$.  So, exactly one of 
$N_{ef}$ and $N_{ef}^{c}$ is a member of $\cal F$.  If it is $N_{ef}$
(resp. $N_{ef}^{c}$), we say that the nonstandard extremities
${\bf e}=[e_{n}]$ and ${\bf f}=[f_{n}]$ are {\em shorted together},
and we write ${\bf e}\asymp {\bf f}$ (resp. ${\bf e}$ and ${\bf f}$
are {\em not shorted together}, and we write ${\bf e}\not\asymp {\bf f}$).
Also we take it that ${\bf e}$
is shorted to itself, i.e., 
${\bf e}\asymp {\bf e}$.  This shorting is an equivalence relation
on the set of nonstandard extremities, whose transitivity is shown by 
\[ \{n\!: e_{n}\asymp f_{n}\}\,\cap\,\{n\!: f_{n}\asymp g_{n}\}\,\subseteq\{n\!: e_{n}\asymp g_{n}\} \]
as usual.  The resulting equivalence classes are the {\em nonstandard
$\omega$-nodes}; typically, they will be denoted by ${\bf x}^{\omega}$. 

Note that each nonstandard $\omega$-node may or may not have a 
nonstandard exceptional element.  In the event that it does have one, say, 
${\bf x}=[e_{n}]$, where the $e_{n}$ are standard nodes 
$x_{n}^{\mu_{n}}$ of natural-number ranks $\mu_{n}$ for almost all $n$,
there are two cases to consider.  In the first case, the ranks 
$\mu_{n}$ assume only finitely many values for almost all $n$.
As was argued in Section 2, there will be exactly one rank $\rho$ 
for which $\{n\!: \mu_{n}=\rho\}\in{\cal F}$.  This allows 
us to identify the rank of $\bf x$ as being $\rho$, and we now write
${\bf x}^{\rho}$ for that nonstandard exceptional element $\bf x$.

The other case arises when, for every $N\in{\cal F}$, the set 
$\{ \mu_{n}\!:n\in N\}$ assumes infinitely many values.
In this case, we can identify the rank $\rho=[\mu_{n}]$ of $\bf x$
as being a hypernatural number that is not a standard number, 
and we may denote $\bf x$ by ${\bf x}^{\rho}$ again.

Here, too, every nonstandard $\omega$-node ${\bf x}^{\omega}$ 
possesses at least one nonstandard $\vec{\omega}$-tip. Also, 
if ${\bf x}^{\omega}$ possesses an exceptional element ${\bf x}^{\rho}$,
that node ${\bf x}^{\rho}$ will be shorted to at least one nonstandard
$\vec{\omega}$-tip, and moreover ${\bf x}^{\rho}$ will not be the 
nonstandard exceptional element of any other nonstandard
$\omega$-node. Furthermore, ${\bf x}^{\omega}$ cannot have two or more 
different nonstandard exceptional 
elements.  These facts, too, follow directly from the 
properties of standard $\omega$-nodes.

Finally, let $^{*}\!X^{\omega}$ denote the set of nonstandard $\omega$-nodes
induced by the originally assumed sequence $\langle G_{n}^{\omega}\rangle$
of standard $\omega$-graphs.  We may now define 
the nonstandard $\omega$-graph 
$^{*}\!G^{\omega}$ induced by $\langle G_{n}^{\omega}\rangle$ to be the 
set\footnote{Since we have assumed that no $G_{n}^{\omega}$ 
possesses any $\vec{\omega}$-nodes, a set 
$^{*}\!X^{\vec{\omega}}$ of nonstandard $\vec{\omega}$-nodes does 
not appear in the set (\ref{4.1}).}
\begin{equation}
^{*}\!G^{\omega}\;=\;\{\,^{*}\!X^{0},\,^{*}\!B,\,^{*}\!X^{1},\ldots,\,^{*}\!X^{\mu},\ldots,\,^{*}\!X^{\omega}\}  \label{4.1}
\end{equation}

\section{Nonstandard Resistive Electrical Networks}

Let $^{*}\!G^{\nu}$ be a nonstandard graph of rank $\nu$, where
$1\leq \nu\leq \omega$, induced by a sequence $\langle G_{n}^{\nu}\rangle$
of standard $\nu$-graphs.  $^{*}\!G^{\nu}$ can be 
converted into a nonstandard, resistive, electrical network
$^{*}\!{\bf N}^{\nu}$ by assigning positive resistances to 
every branch of every $G_{n}^{\nu}$ and voltage sources to some of 
those branches, which we take to be in the Thevenin form (i.e., 
each branch voltage source is connected in series with the branch 
resistance).  This induces hyperreal-valued branch resistances
${\bf r}_{\bf b}$ and hyperreal branch voltage sources ${\bf e}_{\bf b}$ 
in the nonstandard branches $\bf b$ of 
$^{*}\!G^{\nu}$, where again in any branch the hyperreal 
voltage source is connected in series with the hyperreal resistance if that 
source exists there.  This inducement of $^{*}\!{\bf N}^{\nu}$
is exactly the same as that for the nonstandard network 
$^{*}\!{\bf N}^{1}$ of rank 1 explained in \cite[pages 165-166]{gn}.

The question arises as to whether $^{*}\!{\bf N}^{\nu}$
has a hyperreal operating point, that is, a set of hyperreal branch currents 
${\bf i}_{\bf b}$ and a set of hyperreal branch voltages ${\bf v}_{\bf b}$
satisfying in each nonstandard branch $\bf b$ Ohm's law:
\[ {\bf v}_{\bf b}\,=\,{\bf r}_{\bf b}{\bf i}_{\bf b}\,-\,{\bf e}_{\bf b} \]
and Tellegen's equation in a certain way.
The answer is "yes," and it is established in virtually the same way 
as it was established for nonstandard 1-networks $^{*}\!{\bf N}^{1}$
in \cite[Section 8.9]{gn}.  In fact, Theorem 8.9-2 holds word-for-word 
with $^{*}\!{\bf N}^{\nu}$ as it does for $^{*}\!{\bf N}^{1}$.
Now, however, the solution space $^{*}\!{\cal L}$ is based upon
loops of ranks up to $\nu$ instead of up to 1.
Similarly, nonstandard versions of Kirchhoff's current law and 
Kirchhoff's voltage law\footnote{Please see the website 
www.ee.sunysb.edu/\~\,zeman (in particular, the Errata for 
\cite{gn}) for a rather obvious correction to pages 167 and 168 of that book.}
again hold as stated in Theorems 8.9-3 and 8.9-4 respectively, with 
$^{*}\!{\bf N}^{1}$ replaced by $^{*}\!{\bf N}^{\nu}$.
Since all of this is virtually identical to the development given in 
\cite[Section 8.9]{gn}, we will say no more about it, other that to note 
that these nonstandard results for linear networks
can be extended to nonlinear resistive 
networks by exploiting Duffin's theorems noted in the last paragraph
of that Section 8.9 of \cite{gn}.

\section{Nonstandard Graphs and Networks of Still Higher Ranks}

Let us briefly remark that all our results can be extended 
to still higher ranks.  The results for successor-ordinal ranks
can be obtained by mimicking our development for natural-number ranks.
For limit-ordinal ranks, the development mimics that given above 
for rank $\omega$.

\end{document}